\documentclass[a4paper,12pt]{amsart}
\usepackage[arrow,matrix]{xy}
\usepackage{amsmath,amssymb,amscd,amsthm,mathrsfs,hyperref}
\theoremstyle{plain}
\newtheorem{thm}{Theorem}[section]
\newtheorem{cor}[thm]{Corollary}
\newtheorem{lem}[thm]{Lemma}
\newtheorem{prop}[thm]{Proposition}
\newtheorem{defn}[thm]{Definition}
\newtheorem{exa}[thm]{Example}

   \def\op{\oplus} \def\ot{\otimes}
\def\Hom{\operatorname {Hom}}
\def\RHom{\operatorname {RHom}}
\def\Ext{\operatorname {Ext}}

\begin{document}
\title[Hopf algebra actions on DGA]{\bf Hopf algebra actions on differential graded algebras and applications}

\author{Ji-Wei He, Fred Van Oystaeyen and Yinhuo Zhang}
\address{J.-W. He\newline \indent Department of Mathematics, Shaoxing College of Arts and Sciences, Shaoxing Zhejiang 312000,
China\newline \indent Department of Mathematics and Computer
Science, University of Antwerp, Middelheimlaan 1, B-2020 Antwerp,
Belgium} \email{jwhe@usx.edu.cn}
\address{F. Van Oystaeyen\newline\indent Department of Mathematics and Computer
Science, University of Antwerp, Middelheimlaan 1, B-2020 Antwerp,
Belgium} \email{fred.vanoystaeyen@ua.ac.be}
\address{Y. Zhang\newline
\indent Department WNI, University of Hasselt, Universitaire Campus,
3590 Diepenbeek, Belgium} \email{yinhuo.zhang@uhasselt.be}

%\thanks{}
\date{}
\begin{abstract} Let $H$ be a finite dimensional semisimple
Hopf algebra, $A$ a differential graded (dg for short) $H$-module
algebra. Then the smash product algebra $A\#H$ is a dg algebra. For
any dg $A\#H$-module $M$, there is a quasi-isomorphism of dg
algebras: $\RHom_A(M,M)\#H\longrightarrow \RHom_{A\#H}(M\ot H,M\ot
H)$. This result is applied to $d$-Koszul algebras, Calabi-Yau
algebras and
AS-Gorenstein dg algebras.\\[2mm]
{\bf Keywords}: differential graded algebra, smash product, Yoneda algebra\\
{\bf MSC(2000)}: 16E45, 16E40, 16W50
\end{abstract}
\maketitle

\section{Introduction}

Let $G$ be a finite group, $R$ be a $G$-group algebra. For modules
$M$ and $N$ over the skew group algebra $R* G$, there are
``natural'' $G$-actions on the extension groups $\Ext^*_R(M,N)$ and
the Hochschild cohomologies of $R$ (see \cite{M,MMM}). Naturally,
the group actions on extension groups or Hochschild cohomologies can
be generalized to Hopf algebra actions \cite{F,S}. It seems not easy
to find an explicit relations between the extensions of modules over
$R*G$ and that of modules over $R$. Let $R$ be a positively graded
algebra and $G$ be a finite group of grading preserving
automorphisms of $R$. In \cite{M}, the author established an
interesting link between the Yoneda algebra of the trivial module
over $R*G$ and that of the trivial module over $R$. In \cite{MMM},
the authors discussed the Hochschild cohomology algebra of skew
group algebras. In both papers, the results were proved through
discussing the structures of extension groups directly.

Since the extension groups of modules and Hochschild cohomologies
can be computed through differential graded (dg, for short) modules
and dg algebras, we try to deal with the group actions on extensions
under the framework of dg settings, and we find that this is an
efficient way to do these things. So we first need to discuss
properties of the group (more generally, Hopf algebra) actions on dg
algebras. Let $H$ be a finite dimensional semisimple Hopf algebra, and $A$ be a dg
$H$-module algebra. Our main result (Theorem \ref{thm1}) says that
there is a quasi-isomorphism of dg algebras:
$\RHom_A(M,M)\#H\longrightarrow \RHom_{A\#H}(M\ot H,M\ot H)$ for any
dg $A\#H$-module $M$. This quasi-isomorphism yields the isomorphism
in \cite[Theorem 10]{M} and generalizes \cite[Theorem 2.3]{VZ} to
the level of derived functors. We apply the main result to
$d$-Koszul algebras, Calabi-Yau algebras and AS-Gorenstein algebras.
We show that the smash product algebra of a $d$-Koszul algebra and
finite dimensional semisimple Hopf algebra is also a $d$-Koszul
algebra, and the Galois covering algebras of a $d$-Koszul Calabi-Yau
algebras are also Calabi-Yau.

Throughout, $k$ is an algebraically closed field of characteristic
zero and all algebras are $k$-algebras; unadorned $\ot$ means
$\ot_k$ and Hom means Hom$_k$. By a dg algebra we mean a cochain dg
algebra, that is, a graded algebra $A=\bigoplus_{n\in \mathbb{Z}}A^n$
with a differential $d$ of degree 1, such that for all homogeneous
elements $a,b\in A$ we have $d(ab)=d(a)b+(-1)^{|a|}ad(b)$ where
$|a|$ denotes the degree of $a$. An associative algebra $R$ may be
regarded as a dg algebra concentrated in degree zero with zero
differentials, and then a complex of $R$-modules may be regarded as
a dg $R$-module. A (left) dg $A$-module is a graded $A$-module
$M=\op_{n\in\mathbb{Z}}M^n$ with a differential $d$ of degree 1 such
that for all homogeneous elements $a\in A$ and $m\in M$ we have
$d(am)=d(a)m+(-1)^{|a|}ad(m)$. Similarly, we have right dg modules.
In this paper, by a dg module we always mean a left dg module.

Let $A$ be a dg algebra, $M$ and $N$ dg $A$-modules. We write
$\Hom_A(M,N)=\op_{i\in\mathbb{Z}}\Hom_A^i(M,N)$, where $\Hom^i_A(M,N)$
is the set of all graded $A$-module maps of degree $i$. Then
$\Hom_A(M,N)$ is a complex with the canonical differential $d$ which
acts on a homogeneous element $f\in\Hom_A(M,N)$ by $d(f)=d_N\circ
f-(-1)^{|f|}f\circ d_M$. A dg $A$-module $P$ is said to be {\it
homotopically projective} (or K-projective) if $\Hom_A(P,-)$
preserves the quasi-trivial dg modules, and a dg $A$-module $I$ is
said to be {\it homotopically injective} (or K-injective) if
$\Hom_A(-,I)$ preserves the quasi-trivial dg modules (see
\cite[Chapter 8]{KZ}).

For more properties of dg algebras and modules, we refer to the
references \cite{AFH,FHT,K1,KM}.

\section{Hopf algebra actions on dg algebras}

Let $A$ be a dg algebra, $H$ a Hopf algebra. We call $A$ a {\it dg
$H$-module algebra} if
\begin{itemize}
  \item [(i)] $A$ is a graded $H$-module algebra, that is; for $a\in A^i$,
$b\in A^j$ and $h\in H$, $h\cdot a\in A^i$ and
$h\cdot(ab)=(h_{(1)}\cdot a)(h_{(2)}\cdot b)$,  and
\item[(ii)] the differential $d$ of $A$ is compatible with the
$H$-module action, that is, $d(h\cdot a)=h\cdot d(a)$ .
\end{itemize}

If $A$ is a left dg $H$-module algebra, then the cohomology algebra
$H(A)$ is a graded $H$-module algebra, and the smash product $A\#H$
is also a dg algebra with the differential $\delta=d\ot id$.

\begin{prop} If $A$ is a dg $H$-module algebra, then $H(A\#H)\cong H(A)\#H$ as graded algebras. \end{prop}

Since $A$ is a dg subalgebra of $A\#H$, any dg $A\#H$-modules $M$
and $N$ can be viewed as dg $A$-modules. Also $M$ and $N$ are
complexes of $H$-modules. If the antipode $S$ of $H$ is bijective,
there is a natural $H$-module structure on the complex
$\Hom_A(M,N)$. Explicitly, for $f\in\Hom_A(M,N)$, and $h\in H$, the
$H$-module action is defined by
\begin{equation}\label{eq2}
(h\rightharpoonup f)(m)=h_{(2)}f(S^{-1}h_{(1)}m),\qquad \text{for}\
m\in M,\ h\in H.
\end{equation}

Given an $H$-module $X$, we write $X^H=\{x\in X|hx=\varepsilon(h)x,\
\text{for all}\ h\in H\}$ for the invariant submodule. If $M$ is a
complex of $H$-modules, then $M^H$ is a subcomplex of $M$.

\begin{lem}\label{prop1}Let $H$ be a Hopf algebra with a bijective
antipode.

{\rm(i)} If $A$ is a dg $H$-module algebra, $M$ and $N$ are dg
$A\#H$-modules, then $$\Hom_{A\#H}(M,N)\cong\Hom_A(M,N)^H$$ as
complexes of vector spaces

{\rm(ii)} If $H$ is semisimple, then $(\ )^H$ preserves exact
sequences.
\end{lem}
\proof The assertion (i) follows directly from the definition. For
the assertion (ii), just observe that $(\ )^H\cong\Hom_H(k,-)$ is
exact. \qed

Let $M$ be a dg $A$-module. Then $\Hom_A(M,M)$ is a dg algebra. For
convenience, the multiplication of $\Hom_A(M,M)$ is defined as
follows: for homogeneous elements $f,g\in\Hom_A(M,M)$,
$f*g=(-1)^{|f||g|}g\circ f$.

From now on, $H$ will always be a Hopf algebra with a bijective
antipode.

\begin{lem}\label{prop2} Let $A$ be a dg $H$-module algebra, $M$ a dg $A\#H$-module. Then $B=\Hom_A(M,M)$ is
a dg $H$-module algebra.
\end{lem}
\proof Straightforward. \qed

Let $M$ be a dg $A\#H$-module, $W$ a dg $H$-modules. Then $M\ot W$
is a dg $A\#H$-module through the action defined by
\begin{equation}\label{eq1}
    (a\#h)(m\ot w)=(a\#h_{(1)})m\ot h_{(2)}w.
\end{equation}

\begin{lem} Let $M$ and $N$ be dg $A\#H$-modules, $W$ a dg $H$-modules. Then the Hom-Tensor adjoint isomorphism
$$\varphi:\Hom(W,\Hom_A(M,N))\longrightarrow \Hom_A(M\ot W,N)$$ is an $H$-module morphism. \end{lem}
\proof For homogeneous elements $f\in\Hom(W,\Hom_A(M,N))$ and $m\in
M$, $w\in W$, recall that $\varphi(f)(m\ot w)=(-1)^{|m||w|}f(w)(m)$.
Hence for $h\in H$, we have
$$\begin{array}{cll}
\varphi(h\rightharpoonup f)(m\ot w) & = & (-1)^{|m||w|}(h\rightharpoonup f)(w)(m) \\
                                    & = & (-1)^{|m||w|}[h_{(1)}\rightharpoonup(f(S^{-1}(h_{(1)})w))](m) \\
                                    & = & (-1)^{|m||w|}h_{(3)}(f(S^{-1}(h_{(1)})w))(S^{-1}(h_{(2)})m), \\
                                                  \end{array}
$$
$$\begin{array}{ccl}
    (h\rightharpoonup\varphi(f))(m\ot w) & = & h_{(2)}\varphi(f)(S^{-1}(h_{(1)})(m\ot w)) \\
     & = & h_{(3)}\varphi(f)(S^{-1}(h_{(2)})m\ot S^{-1}(h_{(1)})w)) \\
     & = & (-1)^{|m||w|}h_{(3)}(f(S^{-1}(h_{(1)})w))(S^{-1}(h_{(2)})m).
  \end{array}
$$
Therefore $\varphi$ is an $H$-module morphism. \qed

\begin{prop}\label{prop3} Let $H$ be a finite dimensional semisimple Hopf algebra,
$A$ a dg $H$-module algebra. Then a dg $A\#H$-module $P$ is
K-projective if and only if $P$ is K-projective as a dg $A$-module.
\end{prop}
\proof Assume that $P$ is a K-projective $A$-module. We only need to
show that the functor $\Hom_{A\#H}(P,-)$ preserves the quasi-trivial
dg modules. By Lemma \ref{prop1}, we have
$$\Hom_{A\#H}(P,-)=(\ )^H\circ\Hom_A(P,-).$$ Since $H$ is
semisimple, by Lemma \ref{prop2} the functor $(\ )^H$ preserves
exact sequences. Hence $(\ )^H\circ\Hom_A(P,-)$ preserves
quasi-trivial dg modules.

Conversely, suppose that $P$ is a K-projective $A\#H$-module. Since
$P$ is homotopically equivalent to a semifree dg $A\#H$-module (see
\cite{AFH,FHT}), we may assume that $P$ is semifree. Let
$$0\subseteq P(0)\subseteq P(1)\subseteq\cdots\subseteq
P(n)\subseteq P(n+1)\subseteq\cdots$$ be a semifree filtration of
the dg $A\#H$-module $P$. Then $P(n+1)/P(n)$ is a free dg
$A\#H$-module (i.e., it is a direct sum of shifts of $A\#H$).
Moreover, $A\#H$ is also a free dg $A$-module. Hence, the above
filtration is also a semifree filtration of the dg $A$-module $P$.
So $P$ is K-projective as a dg $A$-module. \qed

\begin{prop} Let $A$ and $H$ be as above. A dg $A\#H$-module $I$ is
K-injective if and only if it is K-injective as a dg $A$-module.
\end{prop}
\proof Assume that $I$ is a K-injective dg $A\#H$-module. We have to
show that the functor $\Hom_A(-,I)$ preserves quasi-trivial dg
$A$-modules. As a dg $A$-module, we have
$I\cong\Hom_{A\#H}(A\#H,I)$. Hence we get
$$\Hom_A(-,I)\cong\Hom_A(-,\Hom_{A\#H}(A\#H,I))\cong\Hom_{A\#H}(A\#H\ot_A-,I).$$
Since $A\#H$ is a free right dg $A$-module, $A\#H\ot_A-$ preserves
quasi-trivial dg $A$-modules. Therefore $\Hom_A(-,I)$ preserves the
quasi-trivial dg modules.

The other direction follows from the following isomorphism:
$$\Hom_{A\#H}(-,I)\cong(\ )^H\circ\Hom_A(-,I). \qed$$

\begin{lem} \label{lem1} Let $P$ and $Q$ be dg $A\#H$-modules. If $H$ is finite dimensional, then there is a natural
isomorphism of complexes of vector spaces: $$\Hom_{A\#H}(P\ot H,Q\ot
H)\cong\Hom_A(P,Q)\ot H.$$
\end{lem}
\proof By Lemma \ref{prop1}, we have $$\Hom_{A\#H}(P\ot H,Q\ot
H)\cong\Hom_{A}(P\ot H,Q\ot H)^H.$$ On the other hand, we have an
isomorphism of complexes of $H$-modules
$$\Hom_A(P\ot H,Q\ot H)\cong\Hom(H,\Hom_A(P,Q\ot H)).$$  Hence
$$\begin{array}{ccl}
\Hom_{A\#H}(P\ot H,Q\ot H) & \cong & \Hom(H,\Hom_A(P,Q\ot H))^H \\
 & \cong & \Hom_H(H,\Hom_A(P,Q\ot H)) \\
  & \cong & \Hom_A(P,Q\ot H)\\
  & \cong & \Hom_A(P,Q)\ot H. \qed
       \end{array}
$$

We may write out explicitly the isomorphism in the lemma above as
$$\theta: \Hom_A(P,Q)\ot H\longrightarrow \Hom_{A\#H}(P\ot H,Q\ot
H),$$ acting on elements as $$\theta(f\ot h)(p\ot
g)=g_{(2)}f(S^{-1}(g_{(1)})p)\ot g_{(3)}h,$$ where $f\in
\Hom_A(P,Q)$, $p\in P$ and $g,h\in H$.

Let $M$ be a dg $A\#H$-module. Let $P$ be a K-projective resolution
of the dg $A\#H$-module $M$. From Proposition \ref{prop3}, it
follows that $P$ is also K-projective as a dg $A$-module. Then
$\RHom_A(M,M)=\Hom_A(P,P)$, and hence is a dg algebra. By Lemma
\ref{prop2}, $\RHom_A(M,M)$ is a dg $H$-module algebra. Of course
the dg $H$-module algebra structure of $\RHom_A(M,M)$ depends on the
choice of the K-projective resolution of $M$. However the dg algebra
structures on $\RHom_A(M,M)$ induced from different K-projective
resolutions are quasi-isomorphic to each other as
$A_\infty$-algebras. This does not matter since such dg algebras
have the same homological properties. Also the $H$-module structures
are compatible with the associated quasi-isomorphisms.

\begin{thm} \label{thm1} Let $H$ be a finite dimensional semisimple Hopf
algebra, $A$ a dg $H$-module algebra. If $M$ is a dg $A\#H$-module,
then there is a quasi-isomorphism of dg algebras:
$$\RHom_A(M,M)\#H\longrightarrow \RHom_{A\#H}(M\ot H,M\ot H).$$
\end{thm}
\proof Let $P$ be a K-projective resolution of the dg $A\#H$-module
$M$. Then dg $A\#H$-module $P\ot H$ is quasi-isomorphic to the dg
$A\#H$-module $M\ot H$. By Proposition \ref{prop3}, $P\ot H$ is a
K-projective dg $A\#H$-module. Hence we have
$\RHom_{A\#H}(M,M)=\Hom_{A\#H}(P\ot H,P\ot H)$ and
$\RHom_A(M,M)=\Hom_A(P,P)$. By Lemma \ref{lem1}, we have a
quasi-isomorphism of complexes:
$$\theta: \Hom_A(P,P)\ot H\longrightarrow \Hom_{A\#H}(P\ot H,P\ot
H),$$ and for $f\in \Hom_A(P,P)$, $h,g\in H$ and $p\in P$,
$$\theta(f\ot h)(p\ot g)=g_{(2)}f(S^{-1}(g_{(1)})p)\ot g_{(3)}h.$$ We claim
that $\theta$ is a morphism of dg algebras.

For homogeneous elements $f,f'\in\Hom_A(P,P)$, $h,h',g\in H$ and
$p\in P$, we have $$\begin{array}{ccl}
                      &&\theta((f'\#h')(f\#h))(p\ot g)\\ & = & \theta(f'*(h'_{(1)}\rightharpoonup f)\# h'_{(2)}h)(p\ot g) \\
                       & = & g_{(2)}[f'*(h'_{(1)}\rightharpoonup f)](S^{-1}(g_{(1)})p)\ot g_{(3)}h'_{(2)}h \\
                       & = & (-1)^{|f||g|}g_{(2)}(h'_{(1)}\rightharpoonup f)\left(f'(S^{-1}(g_{(1)})p\right)\ot
                       g_{(3)}h'_{(2)}h,
                    \end{array}
$$ and
$$\begin{array}{ccl}
    &&[\theta(f'\#h')*\theta(f\#h)](p\ot g)\\ & = & (-1)^{|f||g|}\theta(f\#h)\circ\theta(f'\#h')(p\ot g) \\
     & = & (-1)^{|f||g|}\theta(f\#h)(g_{(2)}f'(S^{-1}(g_{(1)})p)\ot g_{(3)}h') \\
     & = & (-1)^{|f||g|}
     g_{(4)}h'_{(2)}f[S^{-1}(g_{(3)}h'_{(1)})g_{(2)}f'(S^{-1}(g_{(1)})p)]\ot
     g_{(5)}h'_{(3)}h\\
     & = &(-1)^{|f||g|}g_{(2)}(h'_{(1)}\rightharpoonup f)\left(f'(S^{-1}(g_{(1)})p\right)\ot
                       g_{(3)}h'_{(2)}h.
  \end{array}
$$
Hence $\theta$ is compatible with the multiplications, and by Lemma
\ref{lem1} it is a quasi-isomorphism. \qed

Let $M$ and $N$ be dg $A\#H$-modules. Since $\RHom_A(M,N)$ is a
complex of $H$-modules, the extension group
$\Ext^*_A(M,N)=\op_{i\in\mathbb{Z}}\Ext^i_A(M,N)$ is a graded
$H$-module, and $\Ext^*_A(M,M)$ is a graded $H$-module algebra.

\begin{cor} \label{cor1} Let $A$ and $H$ be as above, $M$ and $N$ be dg $A\#H$-modules.
\begin{itemize}
  \item [(i)] $\Ext^*_{A\#H}(M,N)\cong\Ext_A^*(M,N)^H$;
  \item [(ii)] $\Ext^*_{A\#H}(M,M)\cong\Ext_A^*(M,M)^H$ as graded
  algebras;
  \item [(iii)] $\Ext^*_{A\#H}(M\ot H,M\ot H)\cong\Ext_A^*(M,M)\#H$ as graded algebras.
\end{itemize}
\end{cor}
\proof (i) Let $P$ and $Q$ be K-projective resolutions of the dg
$A\#H$-modules $M$ and $N$ respectively. Then
$\RHom_{A\#H}(M,N)=\Hom_{A\#H}(P,Q)\cong\Hom_A(P,Q)^H$. Hence
$$\begin{array}{ccl}
    \Ext^*_{A\#H}(M,N) & = & \text{H}^*(\RHom_{A\#H}(M,N))\cong \text{H}^*(\Hom_A(P,Q)^H) \\
     & \cong &(\text{H}^*\Hom_A(P,Q))^H\cong\Ext_A^*(M,N)^H.
  \end{array}
$$

The assertion (ii) is directly from (i). Then assertion (iii) is a
direct consequence of Theorem \ref{thm1}. \qed

Group actions on extension groups have been discussed by several
authors (see \cite{M,MMM}) by the use of the traditional homological
tools. In Corollary \ref{cor1}, when $H=kG$ is a group algebra of
finite group $G$, then assertion (ii) becomes Proposition 2.6 in
\cite{MMM}, and assertion (iii) becomes Theorem 10 in \cite{M}.
Moreover, assertion (iii) is a generalization of \cite[Theorem
2.3]{VZ} to the level of derived functors.

\section{Applications}

Throughout this section, $H$ is a finite dimensional semisimple Hopf
algebra.

\subsection{Hopf algebra actions on $d$-Koszul algebras}

Let $R=\op_{n\ge0}R_n$ be a positively graded algebra such that
$R_0$ is semisimple, dim$R_i<\infty$ for all $i\ge0$ and
$R_iR_j=R_{i+j}$. Recall that $R$ is called a homogeneous algebra if
$R\cong T_{R_0}(R_1)/I$, where $I$ is an ideal generated by elements
in $\underbrace{R_1\ot_{R_0}\cdots\ot_{R_0}R_1}_{d\ factors}$
($d\ge2$). $R$ is called a {\it connected} graded algebra if
$R_0\cong k$. A homogeneous algebra $R$ is called a {\it $d$-Koszul}
algebra if the trivial module $R_0$ has a graded projective
resolution
$$\cdots \longrightarrow P^{-n}\longrightarrow
P^{-n+1}\longrightarrow\cdots\longrightarrow P^0\longrightarrow
R_0\longrightarrow0,$$ such that the graded module $P^{-n}$ is
generated in degree $\frac{n}{2}d$ if $n$ is even and
$\frac{n-1}{2}d+1$ if $n$ is odd for all $n\ge0$. When $d=2$, then a
$d$-Koszul algebra is usually called a {\it Koszul} algebra which
was introduced by Priddy in \cite{P}. The concept of $d$-Koszul
algebra was introduced by Berger in \cite{B}, where a $d$-Koszul
algebra is called a generalized Koszul algebra. Many interesting
algebras are proved to be $d$-Koszul algebras. For example,
3-dimensional graded Calabi-Yau algebras \cite{Bo} are $d$-Koszul
algebras.

For simplicity, write $E^i(R)=\Ext_R^i(R_0,R_0)$ and
$E(R)=\op_{i\ge0}\Ext_R^i(R_0,R_0)$. Endowed with the Yoneda
product, $E(R)$ is a graded algebra. We call $E(R)$ sometimes the
{\it Yoneda Ext-algebra} of $R$. For the $d$-Koszul algebras, we
have the following properties.

\begin{thm}\label{thm2}\cite{BM,BGS,GMMZ,HL} Let $R$ be as above.
\begin{itemize}
  \item [(i)] $R$ is a Koszul algebra if and only if $E(R)$ is generated
by $E^0(R)$ and $E^1(R)$.
  \item [(ii)] If $R$ is a homogeneous algebra, then $R$ is a $d$-Koszul
algebra ($d\ge3$) if and only if $E(R)$ is generated by $E^0(R)$,
$E^1(R)$ and $E^2(R)$.
\end{itemize}\end{thm}

Applying the main result from the last section, we are able to show
that the $d$-Koszulness of a graded algebra can be lifted to a smash
product of the graded algebra.

\begin{thm}\label{thm3} Let $R$ be a homogeneous algebra. Assume that there is an $H$-action on $R$
so that $R$ is a graded $H$-module algebra. Then $R\#H$ is a
$d$-Koszul algebra if and only if $R$ is a $d$-Koszul algebra.
Moreover, $E(R\#H)\cong E(R)\#H$.
\end{thm}
\proof Since $R$ is an $H$-module algebra, $R$ is a graded
$R\#H$-module. In particular, $R_0$ is an $R\#H$-module. Let
$$P^\bullet:=\qquad\cdots \longrightarrow P^{-n}\longrightarrow
P^{-n+1}\longrightarrow\cdots\longrightarrow P^0\longrightarrow
R_0\longrightarrow0$$ be a graded projective resolution of the
$R\#H$-module $R_0$. Write $B=R\#H$. Then $B_0=R_0\#H$. As a left
graded $B$-module, $B_0\cong R_0\ot H$, where the left $B$-module
structure of $R_0\ot H$ is defined by the equation (\ref{eq1}).
Therefore
$$\cdots \longrightarrow P^{-n}\ot H\longrightarrow P^{-n+1}\ot
H\longrightarrow\cdots\longrightarrow P^0\ot H\longrightarrow
B_0\longrightarrow0$$ is a graded projective resolution of the
$B$-module $B_0$. We have the following isomorphisms of graded
algebras
$$\begin{array}{ccl}
    E(R\#H) & = & \op_{i\ge0}H^i\Hom_{R\#H}(P^\bullet\ot H,P^\bullet\ot H) \\
     & \cong & \op_{i\ge0}H^i(\Hom_R(P^\bullet,P^\bullet)\#H)\\
     & \cong & E(R)\#H.
  \end{array}
$$ It is clear that $E(R\#H)$ is generated by $E^0(R\#H)$, $E^1(R\#H)$ and
$E^2(R\#H)$ if and only if $E(R)$ is generated by $E^0(R)$, $E^1(R)$
and $E^2(R)$. Now the proof follows directly from Theorem
\ref{thm2}. \qed

If $R$ is Koszul and $H$ is the group algebra of a finite group $G$,
then the theorem above implies \cite[Theorem 14]{M}.

\subsection{Calabi-Yau algebras}

Let $R$ be a positively graded algebra, $R^e=R\ot R^{op}$ be the
enveloping algebra. $R$ is called a {\it graded Calabi-Yau algebra}
of dimension $p$ (in the sense of Ginzburg \cite{G}) if (i) $R$ is
homologically smooth, that is, as an $R^e$-module $R$ has a
projective resolution of finite length given by finitely generated
modules; (ii) there is a graded $R$-$R$-bimodule isomorphism
\cite{BT}
\begin{eqnarray}\label{eq3}
  % \nonumber to remove numbering (before each equation)
     \Ext^i_{R^e}(R,R^e)&\cong& \left\{
                                                       \begin{array}{ll}
                                                        0 , & i\neq p, \\
                                                         R(l), &
                                                         i=p,
                                                       \end{array}
                                                     \right.
  \end{eqnarray}
where $l$ is an integer, and $R(l)$ is the shift of $R$.

Let $E$ be a finite dimensional graded algebra. We say that $E$ is
{\it graded symmetric} if there is an integer $n$ and a homogeneous
nondegenerate bilinear form $\langle-,-\rangle:E\times
E\longrightarrow k(n)$ such that $\langle xy,z\rangle=\langle
x,yz\rangle$ and $\langle x,y\rangle=(-1)^{|x||y|}\langle
y,x\rangle$ for all homogeneous elements $x,y,z\in E$.

\begin{prop} \label{prop4} Let $Q$ be a finite quiver. A $d$-Koszul algebra
$R=kQ/I$ is a Calabi-Yau algebra if and only if $E(R)$ is a graded
symmetric algebra.
\end{prop}
\proof Assume that $R$ is a Calabi-Yau algebra of dimension $p$. By
\cite[Lemma 4.1]{K2} the triangulated category $D^b(R)$ is a
Calabi-Yau category, where $D^b(R)$ is the triangulated subcategory
of the derived category of $R$ consisting of complexes whose
cohomology has finite total dimension. Then the Yoneda Ext-algebra
$E(R)=\op_{n\ge0}\Ext_R^n(R_0,R_0)=\op_{0\leq n\leq
p}\Hom_{D^b(R)}(R_0,R_0[n])$ is graded symmetric (see \cite[Sect.
2.6]{K2}, or the appendix of \cite{Bo}).

Conversely, assume $E(R)$ is graded symmetric. Suppose that the
global dimension of $R$ is $p$. By \cite[Theorem 1.2]{BM} or
\cite[Theorem 12.5]{LPWZ}, $R$ is an AS-Gorenstein algebra. By
\cite[Proposition 4.5]{BT} $R$ is a Calabi-Yau algebra of dimension
$p$ if and only if $\varepsilon^{p+1}\circ \phi=id$, where
$\varepsilon$ is the isomorphism of $R$ defined by
$\varepsilon(r)=(-1)^{|r|}r$ for a homogeneous element $r\in R$, and
$\phi$ is the isomorphism of $R$ such that $\phi|_{R_1}$ is the dual
map of the restriction map of the Nakayama automorphism of $E(R)$ to
$E^1(R)$. If $d\ge3$, then $\text{gldim}(R)=p$ must be odd. In this
case, $E(R)$ is exactly a symmetric algebra. Hence the Nakayama
automorphism of $E(R)$ is the identity. Therefore $\phi=id$. Since
$p$ is odd, $\varepsilon^{p+1}=id$. Hence $\varepsilon^{p+1}\circ
\phi=id$. That is, $R$ is a Calabi-Yau algebra. If $d=2$, then
$\text{gldim}(R)=p$ can be any positive integer. Thus if $p$ is odd,
the proof is the same as above. However, if $p$ is even, then the
Nakayama automorphism $\nu$ of $E(R)$ satisfies $\nu(x)=(-1)^{|x|}x$
for homogeneous elements $x\in E(R)$. Hence $\phi=\varepsilon$ and
$\varepsilon^{p+1}\circ \phi=id$. That is, $R$ is a Calabi-Yau
algebra. \qed

Let $G$ be a finite group. Suppose that $R$ is an $\mathbb{N}\times
G$-graded algebra such that $R_{0,e}=k$ and $R_{0,g}=0$ for $g\neq
e$, and $R_{i}=\op_{g\in G}R_{i,g}$ is finite dimensional for all
$i\ge0$. Let $M$ and $N$ be finite generated $\mathbb{N}\times G$-graded
$R$-modules. Then $\Hom_R(M,N)$ is an $\mathbb{N}\times G$-graded vector
space. On the other hand, the $\mathbb{N}\times G$-graded algebra $R$
has a natural $kG^*$-module structure so that it is an ($\mathbb{N}$-)graded $kG^*$-module algebra. Similarly, the graded $\mathbb{N}\times G$-graded $R$-modules $M$ and $N$ can be regarded as graded
$R\#kG^*$-modules. Hence $\Hom_R(M,N)$ is a graded $kG^*$-module
with the module structure given by the equation (\ref{eq2}).
However, the $\mathbb{N}\times G$-grading on $\Hom_R(M,N)$ also induces
a graded $kG^*$-module structure. It is not hard to check that the
two $kG^*$-modules described as above coincide. Similarly, if
$P^\bullet$ and $Q^\bullet$ are complexes of graded $\mathbb{N}\times
G$-graded $R$-modules, then $\Hom_R(P^\bullet,Q^\bullet)$ is a
complex of $\mathbb{N} \times G$-graded vector spaces.

We say that an $\mathbb{N}\times G$-graded algebra $R$ is a {\it
Calabi-Yau algebra of dimension $p$} if the isomorphism in
(\ref{eq3}) also respects the $G$-grading. That is, the isomorphism
$\Ext^p_{R^e}(R,R^e)\cong R(l)$ is also an isomorphism of $G$-graded
spaces.

\begin{cor}\label{cor2} Let $R$ be as above. Then $R$ is a $d$-Koszul Calabi-Yau algebra of dimension $p$, then so is $R\#kG^*$.
\end{cor}
\proof If $R$ is a $d$-Koszul algebra, by Theorem \ref{thm3}
$R\#kG^*$ is a $d$-Koszul algebra. By Theorem \ref{thm1},
$\RHom_{R\#kG^*}(k\ot kG^*,k\ot kG^*)\cong \RHom_R(k,k)\#kG^*$.
Since $R$ is a Calabi-Yau algebra of dimension $p$, the global
dimension of $R$ is $p$. Let
$$P^\bullet:=\qquad 0\longrightarrow
P^{-p}\longrightarrow\cdots\longrightarrow P^{-1}\longrightarrow
P^0\longrightarrow k\longrightarrow0$$ be a minimal $\mathbb{N} \times
G$-graded projective resolution of the trivial module $k$. By the
Koszulness of $R$, the projective module $P^{i}$ is finitely
generated for all $i\leq0$. Now
$$\RHom_{R}(k,k)\#kG^*=\Hom_R(P^\bullet,P^\bullet)\#kG^*,$$ and hence
$$E(R\#kG^*)\cong H^*(\Hom_R(P^\bullet,P^\bullet))\#kG^*\cong
E(R)\#kG^*.$$ Since the Yoneda Ext-algebra $E(R)$ is the cohomology
algebra of the dg algebra $\Hom_R(P^\bullet,P^\bullet)$, which
certainly respects the $\mathbb{N} \times G$-gradings of the projective
modules, $E(R)$ is an $\mathbb{N} \times G$-graded algebra. Since $R$ is
Calabi-Yau and the isomorphism $\Ext^p_{R^e}(R,R^e)\cong R(l)$ is
also an isomorphism of $G$-graded vector spaces, one can check
easily that $E(R)$ is in fact an $\mathbb{N} \times G$-graded symmetric
algebra. Then $E(R)\#kG^*$ must be a graded symmetric algebra. It
follows that $E(R\#kG^*)$ is a graded symmetric algebra. Therefore
$R\#kG^*$ is a Calabi-Yau algebra by Proposition \ref{prop4}. Note
that the global dimensions of $R$ and $R\#kG^*$ are the same. Thus
$R\#kG^*$ is also a Calabi-Yau algebra of dimension $p$ since the
global dimension of $R\#kG^*$ coincides with the Calabi-Yau
dimension. \qed

The above corollary is also a direct consequence of \cite[Theorem
17]{F} since in this case $\Ext^p_R(R,R^e)$ is isomorphic to $R(l)$
both as an $R$-$R$-bimodule and as a left $kG^*$-modules. The method
of \cite{F} is meant to compute the Hochschild cohomology of the
smash product algebra by utilizing the spectral sequence obtained in
\cite{S}. However, when it comes to $d$-Koszul case we only need to
compute the Yoneda Ext-algebra of the smash product $R\#kG^*$.

\begin{exa}{\rm Let $R=k[x,y,z]$ be the polynomial algebra. With the natural grading, $R$ is
a $\mathbb{N} $-graded algebra. It is well known that $R$ is a Koszul
Calabi-Yau algebra of dimension 3. Let $\lambda$ be a primitive
$n$th root of unit, and $G=\{\lambda^i|i\in\mathbb{Z}\}$ be the group
generated by $\lambda$. Set $R_{i,g}=0$ if $g\neq \lambda^i$ and
$R_{i,\lambda^i}=R_i$. Then $R$ is an $\mathbb{N} \times G$-graded
algebra. By Theorem \ref{thm3}, $R\#kG^*$ is a Koszul algebra. In
fact, $R\#kG^*$ is a Galois covering \cite{Gr,MMM} of the polynomial
algebra $k[x,y,z]$. Explicitly, let $Q$ be the quiver

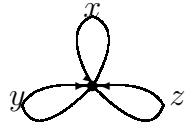
\unitlength=1.2pt
$$\begin{picture}(200,28)(0,0)
\put(100,10){\circle*{3}}
\put(100,34){\makebox(0,0){$x$}}\put(74,4){\makebox{$y$}}\put(124,4){\makebox{$z$}}
\put(102,10){$\vector(-1,0){0}$}\put(98,10){$\vector(1,0){0}$}\put(100.4,11){$\vector(4,-3){0}$}
\qbezier(100,10)(90,28)(100,32)\qbezier(100,10)(110,28)(100,32)
\qbezier(78,4)(82,12)(100,10)\qbezier(78,4)(82,-8)(100,10)
\qbezier(100,10)(118,12)(122,4)\qbezier(100,10)(118,-8)(122,4)
\end{picture}$$ with relations $\rho=\{xy-yx,xz-zx,zy-yz\}$. Let $Q'$
be the following quiver:
$$
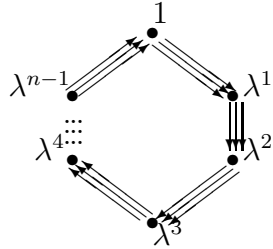
\begin{picture}(200,40)(0,0)
\put(100,30){\circle*{3}}\put(100,33){\makebox{$1$}}
\put(125,10){\circle*{3}}\put(128,10){\makebox{$\lambda^1$}}\put(125,-10){\circle*{3}}\put(128,-10){\makebox{$\lambda^2$}}
\put(100,-30){\circle*{3}}\put(100,-37){\makebox{$\lambda^3$}}
\put(75,-10){\circle*{3}}\put(63,-10){\makebox{$\lambda^4$}}\put(75,10){\circle*{3}}\put(55,10){\makebox{$\lambda^{n-1}$}}
\put(104,28){$\vector(4,-3){22}$}\put(102,27){$\vector(4,-3){22}$}\put(100,26){$\vector(4,-3){22}$}
\put(128,8){$\vector(0,-1){15}$}\put(126,8){$\vector(0,-1){15}$}\put(124,8){$\vector(0,-1){15}$}
\put(123,-12){$\vector(-4,-3){22}$}\put(125,-13){$\vector(-4,-3){22}$}\put(127,-14){$\vector(-4,-3){22}$}
\put(100,-27){$\vector(-4,3){22}$}\put(98,-28){$\vector(-4,3){22}$}\put(96,-29){$\vector(-4,3){22}$}
\put(73,12){$\vector(4,3){22}$}\put(75,11){$\vector(4,3){22}$}\put(77,10){$\vector(4,3){22}$}
\put(72,-5){$\vdots$}\put(74,-5){$\vdots$}\put(76,-5){$\vdots$}
\end{picture}$$
\vspace{8mm}

\noindent For each $0\leq i\leq n-1$, the arrows leaving from the
vertex $\lambda^i$ are labeled as $x_i$, $y_i$ and $z_i$
respectively. The relations of the quiver $Q'$ is
$\rho'=\{x_{i+1}y_i-y_{i+1}x_i,x_{i+1}z_i-z_{i+1}x_i,y_{i+1}z_i-z_{i+1}y_i|0\leq
i\leq n-1\}$, where $x_n=x_0$, $y_n=y_0$ and $z_n=z_0$. Define a map
$F:(Q',\rho')\longrightarrow (Q,\rho)$ of graphs with relations by
sending all the vertices to the unique vertex of $Q$ and sending
arrows $x_i$ to $x$, $y_i$ to $y$ and $z_i$ to $z$ for all $i$. Then
$F$ is a regular covering in the sense of \cite{Gr}. Let
$S=kQ'/(\rho')$ be the quotient algebra of the path algebra $kQ'$ by
modulo the two-side ideal generated by the relations $\rho'$. Then
it is direct to check that $R\# kG^*\cong S$ as $\mathbb{N} $-graded
algebras.

In general, $R\#kG^*$ is not a Calabi-Yau algebra. This is because
the $kG^*$-action on $R$ is not compatible with the Calabi-Yau
property of $R$. In fact, $R$ is Calabi-Yau of dimension 3 as an
$\mathbb{N} \times G$-graded algebra if and only if $\lambda$ is a
 third primitive root of the unit. Now if $G=\{1,\lambda,\lambda^2\}$,
then, by Corollary \ref{cor2}, $R\#kG^*$ is a Calabi-Yau algebra of
dimension 3. The associated quiver $Q'$ is as follows:
$$
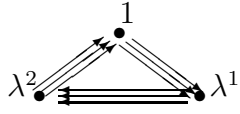
\begin{picture}(200,40)(0,0)
\put(100,30){\circle*{3}}\put(100,33){\makebox{$1$}}
\put(125,10){\circle*{3}}\put(128,10){\makebox{$\lambda^1$}}\put(75,10){\circle*{3}}\put(65,10){\makebox{$\lambda^{2}$}}
\put(104,28){$\vector(4,-3){22}$}\put(102,27){$\vector(4,-3){22}$}\put(100,26){$\vector(4,-3){22}$}
\put(73,12){$\vector(4,3){22}$}\put(75,11){$\vector(4,3){22}$}\put(77,10){$\vector(4,3){22}$}
\put(121,8){$\vector(-1,0){40}$}\put(121,10){$\vector(-1,0){40}$}\put(121,12){$\vector(-1,0){40}$}
\end{picture}$$
By \cite{Bo,Se}, $R\#kG^*$ must be defined by a superpotential
$W'\in kQ'/[KQ',kQ']$ so that $R\#kG^*\cong kQ'/(\partial_{a}W'|a\in
Q'_1)$. In fact, the defining superpotential of $k[x,y,z]$ is
$W=xyz-yxz$, and the defining superpotential of $R\#kG^*$ is the
``lifting'' of $W$, that is, $W'=f_1+f_2+f_3$, where
$f_1=x_2y_1z_0-y_2x_1z_0$, $f_2=z_2x_1y_0-x_2z_1y_0$, and
$f_3=y_2z_1x_0-z_2y_1x_0$.}
\end{exa}

\subsection{Koszul dg algebras and AS-Gorenstein dg algebras}

The concept of Koszul dg algebra was introduced in \cite{HW}. It was
shown that there were some duality properties between a Koszul dg
algebra and its Yoneda Ext-algebra.

\begin{defn}\label{def1}\cite{HW,LPWZ} {\rm Let $A=\op_{n\ge0}A^n$ be a dg algebra such that $A^0$ is semisimple and the differential vanishes
on $A^0$.

(i) $A$ is called a {\it Koszul dg algebra} if $\Ext_A^i(A_0,A_0)=0$
for $i\neq0$.

(ii) $A$ is called a {\it AS-Gorenstein dg algebra} (AS stands for
Artin-Schelter) if there is an integer $n$ such that
$\RHom_{A}(A^0,A)\cong A^0[n]$ as right dg $A$-modules. }
\end{defn}

\begin{prop} Let $A$ be a
dg $H$-module algebra. Then $A$ is a Koszul dg algebra if and only
if $A\#H$ is a Koszul dg algebra.
\end{prop}
\proof Since $A^0$ is semisimple, $A^0\#H$ is also semisimple. By
Corollary \ref{cor1},
$$\Ext^*_{A\#H}(A^0\#H,A^0\#H)\cong\Ext^*_A(A^0,A^0)\#H$$ as graded
algebras. Hence $\Ext^i_{A\#H}(A^0\#H,A^0\#H)=0$ for $i\neq0$ if and
only if $\Ext^i_A(A,A)=0$ for all $i\neq0$. The proof then follows
directly from Definition \ref{def1}. \qed

\begin{prop} Let $A$ be a dg $H$-module algebra such that $A^0=k$. Then $A$ is an AS-Gorenstein dg algebra if and only
if $A\#H$ is an AS-Gorenstein dg algebra.
\end{prop}
\proof By Theorem \ref{thm1},
$\RHom_{A\#H}(A^0\#H,A\#H)\cong\RHom_A(A^0,A)\#H$. Moreover, from
the proof of Theorem \ref{thm1}, one sees that the isomorphism is
compatible with the right $A\#H$-module structures. If $A$ is
AS-Gorenstein, then $\RHom_{A\#H}(A^0\#H,A\#H)\cong A^0\#H[n]$.
Hence $A\#H$ is an AS-Gorenstein algebra. Conversely, if $A\#H$ is
AS-Gorenstein, then $\RHom_A(A^0,A)\#H\cong A^0\#H[n]\cong H[n]$. It
follows that $\Ext_A^i(A^0,A)=0$ for $i\neq n$, and
$\Ext_A^n(A^0,A)$ must be of dimension 1. By suitable truncations of
the right dg $A$-module $\RHom_A(A^0,A)$, we get
$\RHom_A(A^0,A)\cong k[n]$. Therefore $A$ is AS-Gorenstein. \qed

\subsection*{Acknowledgement}  The work is supported by an FWO-grant and NSFC
(No. 10801099).

\bibliography{}

\begin{thebibliography}{99}
\bibitem{AFH} L.L. Avramov, H.-B. Foxby and S. Halperin,
Differential Graded Homological Algebra, preprint.
\bibitem{B} R. Berger, {\it Koszulity for nonquadratic algebras}, J. Algebra {\bf 239} (2001),
705--734.
\bibitem{BM} R. Berger and N. Marconnet, {\it Koszul and Gorenstein properties for homogeneous algebras},
Algebra Represent. Theory {\bf 9} (2006), 67--97.
\bibitem{BT} R. Berger and R. Taillefer, {\it Poincar\'{e}-Birkhoff-Witt deformations of
Calabi-Yau algebras}, J. Noncommut. Geom. {\bf 1} (2007), 241--270.
\bibitem{BGS} A.A. Beilinson, V. Ginzburg and W. Soergel, {\it Koszul duality patterns
               in representation theory}, J. Amer. Math. Soc. {\bf 9} (1996), 473--527.
\bibitem{Bo} R. Bocklandt, {\it Graded Calabi-Yau algebras of dimension
3}, J. Pure Appl. algebra {\bf 212} (2008), 14--32.
\bibitem{F} M. Farinati, {\it Hochschild duality, localization, and smash
products}, J. Algebra {\bf 284} (2005), 415--434.
\bibitem{FHT} Y. F\'{e}lix, S. Halperin and J.-C. Thomas, {Rational Homotopy Theory},
              Grad. Texts Math. {\bf 205}, Springer-Verlag, New York, 2001.
\bibitem{G} V. Ginzburg, {\it Calabi-Yau algebras}, arxiv:math.AG/0612139.
\bibitem{Gr} E.L. Green, {\it Graphs with relations, coverings and group-graded
algebras}, Trans. Amer. Math. Soc. {\bf 279} (1983), 297--310.
\bibitem{GMMZ} E.L. Green, E.N. Marcos, R. Mart\'{\i}nez-Villa and P. Zhang, {\it D-Koszul algebras},
   J. Pure Appl. Algebra {\bf 193} (2004), 141--162.
\bibitem{HL} J.-W. He and D.-M. Lu, {\it Higher Koszul Algebras and A-infinity Algebras}, J. Algebra {\bf 293} (2005),
335--362.
\bibitem{HW} J.-W. He and Q.-S. Wu, {\it Koszul differential graded algebras and BGG correspondence},
             J. Algebra {\bf 320} (2008), 2934--2962.
\bibitem{K1} B. Keller, {\it On differential graded categories}, International Congress of Mathematicians,
{\bf II} (2006), 151--190, Eur. Math. Soc., Z\"{u}rich.
\bibitem{K2} B. Keller, {\it Calabi-Yau triangulated categories},
available at the website:\newline
http://people.math.jussieu.fr/~keller.
\bibitem{KZ} S. K\"{o}nig and A. Zimmermann, Derived Equivalences for
Group Rings, Lect. Notes Math. {\bf 1685}, Spring-Verlag, New York,
1998.
\bibitem{KM} I. K\v r\'{\i}\v z and J.P. May, {\it Operads, algebras, modules and motives},
              Ast\'{e}risque {\bf 233}, 1995.
\bibitem{LPWZ} D.-M. Lu, J. H. Palmieri, Q.-S. Wu and J. J. Zhang,
        {\it $A_\infty$-algebras for ring theorists}, Alg. Colloq. {\bf 11} (2004), 91--128.
\bibitem{M} R. Mart\'{\i}nez-Villa, {\it Skew group algebras and their
Yoneda algebras}, Math. J. Okayama Univ. {\bf 43} (2001), 1--16.
\bibitem{MMM} E.N. Marcos, R. Mart\'{\i}nez-Villa and Ma.I.R. Martins,
{\it Hochschild cohomology of skew group rings and invariants},
Cent. Euro. J. Math. {\bf 2} (2004), 177--190.
\bibitem{P} S. Priddy, {\it Koszul resolutions}, Trans. Amer. Math. Soc. {\bf 152} (1970), 39--60.
\bibitem{Se} E. Segal, {\it The $A_\infty$ deformation theory of a point and the derived categories of local
Calabi-Yaus}, arXiv:math/0702539.
\bibitem{S} D. Stefan, {\it Hochschild cohomology on Hopf Galois
extensions}, J. Pure Appl. Algebra {\bf 103} (1995), 221--233.
\bibitem{VZ} F. Van Oystaeyen and Y. Zhang, {\it H-module endomorphism
rings}, J. Pure Appl. Algebra {\bf 102} (1995), 207--219.


\end{thebibliography}

\end{document}